\def\a{{\alpha}}
\def\b{{\beta}}
\def\d{{\delta}}
\def\e{{\epsilon}}
\def\g{{\gamma}}
\def\k{{\kappa}}
\def\m{{\mu}}
\def\p{{\pi}}

\def\t{{\tau}}
\def\Th{{\Theta}}
\def\w{{\omega}}
\def\W{{\Omega}}

\def\G{{\cal G}}
\def\HH{{\cal H}}
\def\K{{\cal K}}
\def\P{{\cal P}}
\def\Q{{\cal Q}}

\def\bB{{\mathbf B}}
\def\bC{{\mathbf C}}
\def\bE{{\mathbf E}}
\def\bG{{\mathbf G}}

\def\bv{{\mathbf v}}
\def\bZ{{\mathbf Z}}

\def\oa{{\overline\a}}
\def\ob{{\overline\b}}

\def\oB{{\overline B}}
\def\obB{{\overline\bB}}
\def\obv{{\overline\bv}}

\def\Fh{{\hat F}}
\def\Gnp{{\bG_{n,p}}}

\def\inf{{\infty}}
\def\pf{{\Box}}
\def\ra{{\rightarrow}}

\def\l({{\left(}}
\def\r){{\right)}}
\def\'{{\prime}}
\def\lf{{\lfloor}}
\def\rf{{\rfloor}}
\def\({{\Biggl(}}
\def\){{\Biggr)}}
\def\[{{\Biggl[}}
\def\]{{\Biggr]}}
\def\pebb{\atopwithdelims<>}

\def\diam{{\rm diam}}
\def\var{{\rm var}}

\documentclass[12pt]{article}
\newtheorem{thm}{Theorem}

\newtheorem{conj}[thm]{Conjecture}
\newtheorem{cor}[thm]{Corollary}

\newtheorem{lem}[thm]{Lemma}
\newtheorem{prob}[thm]{Problem}

\newtheorem{quest}[thm]{Question}
\newtheorem{res}[thm]{Result}

\usepackage{latexsym,psfig,multicol,amssymb}

\begin{document}

%##############################################################################
%##############################################################################
%
%       TITLE PAGE:
%
\title{Girth, Pebbling, and Grid Thresholds}

\author{Andrzej Czygrinow\\
Department of Mathematics and Statistics\\
Arizona State University\\
Tempe, Arizona 85287-1804\\
email: andrzej@math.la.asu.edu\\
 \\
and \\
 \mbox{}\\
Glenn Hurlbert\thanks{Partially supported by National Security Agency
	grant \#MDA9040210095.}~\\
Department of Mathematics and Statistics\\
Arizona State University\\
Tempe, Arizona 85287-1804\\
email: hurlbert@asu.edu\\
}
\maketitle
\newpage

%#############################################################################
%#############################################################################
%
%       ABSTRACT:
%
\begin{abstract}
In this note we answer a question of Hurlbert about pebbling
in graphs of high girth.
Specifically we show that for every $g$ there is a Class 0 graph of girth
at least $g$.
The proof uses the so-called Erd\H{o}s construction and employs a recent
result proved by Czygrinow, Hurlbert, Kierstead and Trotter.
We also use the Czygrinow et al. result to prove that Graham's pebbling
product conjecture holds for dense graphs.
Finally, we consider a generalization of Graham's conjecture to thresholds 
of graph sequences and find reasonably tight bounds on the pebbling threshold
of the sequence of $d$-dimensional grids, verifying an important
instance the generalization.
\vspace{0.2 in}

\vspace{0.2 in}

\end{abstract}

\newpage

%##############################################################################
%##############################################################################
%
%       BEGINNING OF PAPER:
%
\section{Introduction}\label{Intro}

%##############################################################################
%
%       PEBBLING:
%
\subsection{Pebbling}\label{Pebbling}

A {\it pebbling configuration $\bC$} on a graph $G$ is a distribution of
pebbles to the vertices of $G$.
Given a particular configuration, one is allowed to move the pebbles
about the graph according to this simple rule: if two or more vertices
sit at vertex $v$, then one of them can be moved to a neighbor provided
another is removed from $v$.
Given a specific {\it root} vertex $r$, we say that $\bC$ is $r$-{\it solvable}
if one can more a pebble to $r$ after several pebbling steps, and that
$\bC$ is {\it solvable} if it is $r$-solvable for every $r$.
The pebbling number is the least number $\p=\p(G)$ so that every
configuration of $\p$ pebbles on $G$ is solvable.

The two most obvious pebbling facts are for complete graphs and paths.
The pigeonhole principle implies that $\p(K_n)=n$, and $\p(P_n)=2^{n-1}$
follows by induction or a simple weight function method.
In fact, $\p(G)\ge\min\{n(G),\diam(G)\}$ for every $G$.
Results for trees (a formula based on the minimum path partition of
a tree in \cite{Moe}), $d$-dimensional cubes $Q^d$ (see \cite{Chu}), and
many other graphs with interesting properties are known (see the
survey \cite{Hur}).

An interesting probabilistic version of pebbling was introduced in \cite{CEHK}.
In order to state this variation we introduce some asymptotic notation.
Let $f$ and $g$ be functions of $n$ that tend to infinity.
Denote by $O(g)$ ($o(g)$) the set of functions $f$ for which the ratio
$f/g$ is bounded from above (tends to 0).
Then $g\in\W(f)$ ($g\in\w(f)$) if and only if $f\in O(g)$ ($f\in o(g)$).
We write $f\ll g$ when $f\in o(g)$, $f\sim g$ when $f/g\ra 1$ as $n\ra\inf$,
and set $\Th(f)=\W(f)\cap O(f)$.

Let $\G=(G_i)_{i=1}^\inf$ be a sequence of graphs with strictly increasing
numbers of vertices $N=n(G_i)$.
For a function $t=t(N)$ let $\bC_t$ denote a configuration on $G_i$ that is
chosen uniformly at random from all configurations of $t$ pebbles.
The sequence $\G$ has {\it pebbling threshold} $\t=\t(\G)$ if,
for every $\w\gg 1$,
(1) $\Pr[\bC_t {\rm\ is\ solvable}]\ra 0$ for $t=N/\w$ and
(2) $\Pr[\bC_t {\rm\ is\ solvable}]\ra 1$ for $t=\w N$.

It was proven in \cite{Cla} that the sequence of cliques has threshold
$\t(\K)=\Th(N)$.
Bekmetjev, et al. \cite{BBCH}, showed recently that every graph sequence
has a pebbling threshold.
Bounds on the sequence of paths have undergone several improvements,
the results of which are summarized as follows.
\begin{res}\label{Paths}
The pebbling threshold for the sequence of paths $\P=(P_n)_{n=1}^\inf$
satisfies
$$\t(\P)\in\W\bigg(N2^{c\sqrt{\lg N}}\bigg)\cap O\bigg(N2^{\sqrt{\lg N}}\bigg)$$
for every $c<1/\sqrt{2}$.
\end{res}
The lower bound is found in \cite{BBCH} and the upper bound is found
in \cite{GJSW}.

It is important to draw a distinction between this random pebbling model
and the one in which each of $t$ pebbles independently chooses uniformly
at random a vertex on which to be placed.
In the world of random graphs, the analogs of these two models are
asymptotically equivalent.
However, in the pebbling world, they are vastly different.
For example, in the independent model the pebbling threshold for paths
is at most $N\lg N$ since, with more than that many pebbles,
almost always every vertex already has a pebble on it.

Another important result was proved recently in \cite{CW}.

\begin{res}\label{Cubes}
The pebbling threshold for the sequence of cubes $\Q=(Q^d)_{d=1}^\inf$
satisfies
$$\t(\Q)\in\W(N^{1-\e})\cap O(N)$$
for every $\e>0$.
\end{res}

%##############################################################################
%
%       RESULTS:
%
\subsection{Results}\label{Results}

Pachter et al. \cite{PSV} proved that every graph of diameter two on $N$
vertices has pebbling number either $N$ or $N+1$.
Graphs $G$ with $\p(G)=N(G)$ are called {\it Class 0}, and in \cite{CHH}
a characterization of diameter two Class 0 graphs was found and used to
prove that diameter two graphs with connectivity at least 3 are Class 0.
The authors also conjectured that every graph of fixed diameter and
high enough connectivity was Class 0.
This conjecture was proved by Czygrinow, Hurlbert, Kierstead and Trotter
\cite{CHKT} in the following result.

\begin{res}\label{conn}
Let $d$ be a positive integer and set $k=2^{2d+3}$.
If $G$ is a graph of diameter at most~$d$ and connectivity at least~$k$,
then $G$ is of Class 0.
\end{res}

In this note, we present two applications of this result.
Our first application concerns the following girth problem posed in \cite{Hur}.

\begin{quest}
Does there exist a constant $C$ such that if $G$ is a connected graph
on $n$ vertices with $girth(G) > C$ then $\p(G) > n$?
\end{quest}

Using the so-called Erd\H{o}s construction \cite{Erd}, we answer
the above question in the negative.
Let $g_0(n)$ denote the maximum number $g$ such that there exists
a graph $G$ on $n$ vertices with $girth(G)\geq g$ and $\p(G)=n$.
That is, $g_0(n)$ is the highest girth, as a function of $n$,
among all Class 0 graphs.
It is easy to see that
$$g_0(n) \leq 2\lg{n}$$
(because the cycle on $k$ vertices has pebbling number at least
$2^{\lf k/2\rf}$ --- see \cite{PSV})
and we prove the following lower bound.

\begin{thm}\label{girth}
There exist $n_0$ and $c$ such that, for every $n\geq n_0$,
$$g_0(n) \geq c \sqrt{\lg{n}}\ .$$
\end{thm}
We prove this theorem in Section \ref{GirthProof}.

Our second application concerns the following conjecture of Graham \cite{Chu}.

\begin{conj}\label{Graham}
Every pair of graphs $G$ and $H$ satisfy $\p(G\Box H)\le \p(G)\p(H)$.
\end{conj}

Here, the Cartesian product has vertices $V(G\Box H)=V(G)\times V(H)$ and edges
$E(G\Box H)=\{u\times E(H)\}_{u\in V(G)}\cup\{E(G)\times v\}_{v\in V(H)}$.
A number of theorems have been published in support of this conjecture,
including the recent work of Herscovici \cite{Her} which verifies the case
for all pairs of cycles.
We show the following.

\begin{thm}\label{prod}
Let $G$ and $H$ be connected graphs on $n$ vertices with minimum degrees 
$\d(G)$, $\d(H)$ and let $\d =\min\{\d(G),\d(H)\}$.
If $\d \geq 2^{12n/\d +15}$ then $G \Box H$ is of Class 0.
\end{thm}
In particular, if $\d \gg \frac{n}{\lg{n}}$ then $G \Box H$ is of Class 0. 
We prove this in section \ref{ProdProof}, again using Result \ref{conn}.
As a corollary we obtain that Graham's Conjecture is satisfied for graphs 
with minimum degree $\d \gg \frac{n}{\lg{n}}$.

\begin{cor}\label{grahcor}
Let $G$ and $H$ be such as in Theorem \ref{prod}. 
Then $\p(G\Box H)\le \p(G)\p(H)$.
\end{cor}  

\noindent
{\it Proof.}
We have $\p(G\Box H)=n(G\Box H)=n(G)n(H)\le \p(G)\p(H)$.
\hfill$\pf$
\medskip

Finally, in this paper we also consider the following probabilistic
analog of Graham's Conjecture \ref{Graham},
which we consider a correction of one from \cite{Hur}.

\begin{prob}\label{ThreshProd}
Let $\G=(G_n)_{n=1}^\inf$ and $\HH=(H_n)_{n=1}^\inf$ be two graph sequences.
Define the {\it product sequence} $\G\Box\HH=(G_n\Box H_n)_{n=1}^\inf$.
Find $\t(\G\Box\HH).$
\end{prob}
Let $N(H_n)$, $N(G_n)$ denote the number of vertices of graphs $H_n$ and $G_n$ from Problem \ref{ThreshProd}.
It would be interesting to determine for which sequences 
$\G=(G_n)_{n=1}^\inf$ and $\HH=(H_n)_{n=1}^\inf$, we have
\begin{equation}\label{ThreshProdCon}
\t(\G\Box\HH)\in O\(g\bigg(N(H_n)\bigg)h\bigg(N(G_n)\bigg)\),
\end{equation}
where $g\in\t(\G)$ and $h\in\t(\HH)$.
We call pairs of sequences which satisfy (\ref{ThreshProdCon})
{\it well-behaved}.
One might conjecture that all pairs of sequences are well-behaved,
but we believe counterexamples might exist.

%The formulation in terms of $\a$ and $\b$ has to do with compensating for
%how much of the product graph is accounted for by each term of the product,
%and using the corresponding portion of each threshold
%(notice that we have $N^\oa=N(G_n)$ and $N^\ob=N(H_n)$).
%For example, if $\a=\b$ then the conjectured upper bound reduces to
%$g(\sqrt N)h(\sqrt N)$.

We define the two-dimensional grid $P_n^2=P_n\Box P_n$, and in general the
$d$-dimensional grid $P_n^d=P_n\Box P_n^{d-1}$.
It is easy to show that $P_n^d=P_n^\a\Box P_n^\b$ for all
$\a$ and $\b$ for which $\a+\b=d$.
If we denote $\P^d=(\P_n^d)_{n=1}^\inf$ then we have $\P^d=\P^\a\Box\P^\b$.
Thus, for example, in light of Result \ref{Paths}, the truth of (\ref{ThreshProdCon})
 would imply that
$$\t(\P^2) \in O\( \bigg( \sqrt{N} 2^{ \sqrt{ \lg \sqrt{N} } } \bigg)^2 \)
	= O\( N 2^{ \sqrt{ 2 \lg n } } \) \ .$$
Here we prove the following stronger theorem.

\begin{thm}\label{grid}
Let $\P^d=(\P_n^d)_{n=1}^\inf$ be the sequence of $d$-dimensional grids,
where $P_n^d=(P_n)^d$ is the cartesian product of $d$ paths on $n$ vertices
each, and let $N=n^d$ be the number of vertices of $\P_n^d$.
Then
$$\t(\P^d)\subseteq
	\W\bigg(N2^{c_d(\lg N)^{1/(d+1)}}\bigg)
	\cap O\bigg(N2^{c^\'_d(\lg N)^{1/(d+1)}}\bigg)$$
for all $c_d<2^{-1/2d}$ and $c^\'_d>d+1$.
\end{thm}

This verifies (\ref{ThreshProdCon}) in the case of grids.

\begin{cor}\label{GridCor}
Let $\a$, $\b$ be any pair of integers then for $\G=\P^\a$ and $\HH=\P^\b$, 
(\ref{ThreshProdCon}) holds.
\end{cor}

\noindent
{\it Proof.}
Indeed, if $g\in\t(\G)$ and $h\in\t(\HH)$ then Theorem \ref{grid} says that
\begin{eqnarray*}
g(N^\oa)h(N^\ob)
	& \in & \W\bigg(N^\oa2^{c_\a(\lg N^\oa)^{1/(\a+1)}}
		N^\ob2^{c_\b(\lg N^\ob)^{1/(\b+1)}}\bigg)\\
&&\\
	& \subseteq & \W\bigg(N2^{c(\lg N)^{1/(\g+1)}}\bigg)\\
&&\\
	& \subseteq & \W\bigg(N2^{c(\lg N)^{1/(d/2+1)}}\bigg)\ ,
\end{eqnarray*}
for some $c$, where $\g=\max\{\a,\b\}$ and $d=\a+\b$.
On the other hand, Theorem \ref{grid} also says that
$$\t(\P^{\a+\b})=\t(\P^d)
	\in O \bigg(N2^{c^\'_d(\lg N)^{1/(d+1)}}\bigg)\ ,$$
which is asymptotically smaller.
\hfill$\pf$
\vspace{0.2 in}

We prove Theorem \ref{grid} in Section \ref{GridProof}.

%##############################################################################
%##############################################################################
%
%       PROOFS:
%
\section{Proofs}\label{Proofs}

%##############################################################################
%
%       FIRST THEOREM:
%
\subsection{Proof of Theorem \ref{girth}}\label{GirthProof}

We will need the Chernoff-Hoeffending inequality (see \cite{CHI}).

\begin{res}\label{chern}
Let $X=B(n,p)$ be a binomial random variable with expectation $\m = \bE[X]$.
Then for every $0 < t < \m$,
$$\Pr[|X - \m| > t] < 2 e^{-t^2/3\m}\ .$$
\end{res}

We will also make use of Mader's theorem (see \cite{Mad}), below.

\begin{res}\label{mader}
Every graph having average degree at least $\bar d$
has a subgraph of connectivity at least ${\bar d}/4$.
\end{res}

\noindent
{\it Proof of Theorem \ref{girth}.}
Let $g=g(n) \ll \sqrt{\lg{n}}$, $p =n^{-1 +1/g}$, and consider
the random graph $\bG=\Gnp$ (edges appear in $\bG$ independently
with probability $p$).
Let $X$ denote the number of cycles in $\bG$ of length at most $g-1$.
Then $X$ has expectation
$$\m\ =\ \sum_{i=3}^{g-1}{n \choose i} \frac{(i-1)!}{2}p^i\ 
	<\ g n^{1-1/g}\ .$$
Thus by Markov's inequality (see \cite{CHI}),
\begin{equation}\label{eq1}
\Pr[X > n/4]\ \leq\ 4g/n^{1/g}\ \ra\ 0\ .
\end{equation}

For every vertex $v$ in $\bf{G}$, the random variable $\deg(v)$ has
a binomial distribution and so, by Result \ref{chern},
$$\Pr[|\deg(v)- n^{1/g}| > n^{1/g}/4]\ <\ 2e^{-(n^{1/g})/48}\ .$$
Consequently,
\begin{equation}\label{eq2}
\Pr[ |\deg(v)- n^{1/g}| > n^{1/g}/4\ {\rm for\ some\ } v ]\ \ra\ 0\ .
\end{equation}

Therefore there exists a graph $G$ such that $X \leq n/4$ and,
for every vertex $v$, $|\deg(v)- n^{1/g}| \leq n^{1/g}/4$.
Let $H$ be obtained from $G$ by deleting one vertex from each cycle
of length less than $g$.
Then $|V(H)| \geq 3n/4$ and
$$\sum_{v\in V(H)}\deg(v)\ 
	\geq\ \sum_{v \in V(G)} \deg(v) - \frac{5}{8}n^{1+1/g}\ 
	>\ \frac{1}{8}n^{1+1/g}\ .$$
Thus the average degree of $H$ is at least $n^{1/g}/8$.
By Result \ref{mader}, $H$ contains a subgraph $F$ which is
$n^{1/g}/32$-connected.
Clearly $F$ has girth at least $g$.

Finally, let $\Fh$ be an edge maximal graph on the same vertices as $F$
such that $F$ is a subgraph of $\Fh$ and $\Fh$ has girth at least $g$.
We claim that the diameter of $\Fh$ is at most $g-1$.
Indeed, suppose that there exist $x$ and $y$ such that the shortest path
between $x$ and $y$ has length at least $g$.
Then we can add $xy$ to $\Fh$ to obtain a graph of girth
at least $g$, which contradicts the maximality of $\Fh$.
Therefore, $\Fh$ has girth at least $g$, diameter at most $g$
and is $n^{1/g}/32$-connected.
Since $g \ll \sqrt{\lg{n}}$, we can apply Theorem \ref{conn} to conclude
that $\Fh$ is of Class 0.
\hfill$\pf$
\vspace{0.2 in}

%##############################################################################
%
%       SECOND THEOREM:
%
\subsection{Proof of Theorem \ref{prod}}\label{ProdProof}
Theorem \ref{prod} follows from the following two lemmas and Result \ref{conn}.

\begin{lem}
\label{s_diam}
Let $G$ be a connected graph on $n$ vertices with minimum degree $\d$.
Then the diameter of $G$ is at most $3\frac{n}{\d}+3$.
\end{lem}

\noindent
{\it Proof.}
Fix two vertices $x$, $y$ in $G$ and consider the shortest path 
$x=x_1, \dots, x_k =y$ between $x$ and $y$.
Let $i = \lfloor \frac{k-1}{3} \rfloor$. 
Then $x_1, x_4, x_7 \dots, x_{3i+1}$ must have disjoint neighborhoods, and so
$i (\d+1) \leq n$ which yields $\frac{k-1}{3} -1 < \frac{n}{\d+1}$, so that
$k < \frac{3n}{\d+1}+4 \leq \frac{3n}{\d}+3$.
\hfill$\pf$
\medskip

The next Lemma was proved by Czygrinow and Kierstead \cite{CK}. 
We reproduce the proof here. 

\begin{lem}
\label{s_conn}
The product $G\Box H$ has connectivity
$\k(G \Box H) \geq \min\{\d(G),\d(H)\}$. 
\end{lem}

\noindent
{\it Proof.}
Set $\d=\min\{\d(G),\d(H)\}$.
Let $v_1=(g,h_1), v_2=(g,h_2),\dots , v_{\d}=(g,h_{\d})$, 
$w_1=(g_1,h), w_2=(g_2,h), \dots , w_{\d}=(g_{\d},h)$ 
be distinct vertices in $G\Box H$ that satisfy
\begin{equation}
\label{gcond}
dis_{G}(g_i,g) \leq dist_G(g_{i+1},g)
\end{equation}
and
\begin{equation}
\label{hcond}
dis_{H}(h_i,h) \leq dist_H(h_{i+1},h),
\end{equation}
for $i=1,\dots , \d-1$.
We shall construct vertex-disjoint paths $P_1, \dots , P_{\d}$ 
such that $P_i$ connects $v_i$ with $w_i$. Construct $P_1$ as follows. 
Let $g_1\bar{g}(1)\dots \bar{g}(k)g$ be any shortest path in $G$ 
connecting $g_1$ with $g$ and let $h\bar{h}(1)\dots \bar{h}(l)h_1$ 
be any shortest path in $H$ connecting $h$ with $h_1$.
Then $P_1$ is the path:
$$w_1=(g_1,h)(g_1,\bar{h}(1))\dots (g_1,h_1)(\bar{g}(1),h_1)\dots (g,h_1)=v_1.$$
Delete $v_1$ and $w_1$ and construct $P_2,\dots ,P_{\d}$ inductively. 
We claim that $P_2,\dots , P_{\d}$ are vertex-disjoint with $P_1$. 
Indeed, suppose that $V(P_j)\cap V(P_1)\neq \emptyset$ 
for some $j=2,\dots , \d$. 
There are two similar cases to consider. 
First, suppose that $(g_j, f) \in V(P_j)\cap V(P_1)$. 
Since $g_j \neq g_1$, $f=h_1$ and $g_j=\bar{g}(i)$ for some $i=1, \dots k$. 
Then however
$$dist_G(g_j,g) < dist_G(g_1,g),$$
contradicting (\ref{gcond}). 
Similarly, if $(f,h_j) \in V(P_j)\cap V(P_1)$ then $f=g_1$ and
$h_j=\bar{h}(i)$ for some $i=1, \dots l$ which implies that
$$dist_G(h_j,h) < dist_G(h_1,h),$$ 
contradicting (\ref{hcond}).

By induction, paths $P_1, \dots , P_{\d}$ are vertex-disjoint. 
Now, for any two vertices 
$v=(g, \tilde{h}),w=(\tilde{g},h) \in V(G \Box H)$, let
$v_1= (g,h_1), v_2=(g,h_2),\dots , v_{\d}=(g,h_{\d})$ be 
neighbors of $v$ in $G$-dimension,
$w_1=(g_1,h), w_2=(g_2,h), \dots , w={\d}=(g_{\d},h)$ 
neighbors of $w$ in $H$-dimension ordered accordingly 
to (\ref{gcond}) and (\ref{hcond}).
By the previous argument we can find vertex-disjoint path 
$P_1, \dots, P_{\d}$ 
connecting $v_i$'s with $w_j$'s. 
These paths can be now used to connect 
$v$ with $w$ by $\d$ internally vertex-disjoint paths. 
Indeed, if any of the paths contains $v$ or $w$ then it 
yields a shorter path between $v$ and $w$ which 
is disjoint with other paths. If $v$ is connected with $w$ then there are $\d-1$ 
internally vertex-disjoint path connecting neighbors of $v$ (other than $w$) 
with neighbors of $w$ (other than $v$) and the path $vw$.
Therefore by Menger's Theorem (see \cite{Mad}) 
the connectivity of $G \Box H$ is at least $\d$.

\hfill$\pf$
\medskip

\noindent
{\it Proof of Theorem \ref{prod}.}
By Lemma \ref{s_diam}, the diameter $d$ of $G \Box H$ is at most 
$6\frac{n}{\d}+6$ and by Lemma \ref{s_conn},
the connectivity $k$ of $G\Box H$ is at least $\d$. 
Since $\d \geq 2^{12n/\d +15}$ the assumptions of Result \ref{conn} are 
satisified and so $G \Box H$ is of Class 0.
\hfill$\pf$

%##############################################################################
%
%       THIRD THEOREM:
%
\subsection{Proof of Theorem \ref{grid}}\label{GridProof}

Throughout, we let $N=n^d$.
Also, we define ${a\pebb b}={a+b-1\choose b}$.
Note that $a\pebb b$ is the number of ways to place $b$ unlabeled balls
into $a$ labeled urns.
For our purposes, it equals the number of configurations of $b$ pebbles
on a graph of $a$ vertices.
We will also use the fact that $a\pebb b$ counts the number of points in 
$\bZ^a$ whose coordinates are nonnegative and sum to $b$.

We begin by proving that a configuration with relatively few pebbles
almost always has no vertices having a huge number of pebbles.
For natural numbers $a$ and $b$ define $a^{\underline b}=a!/(a-b)!$.

\begin{lem}\label{flat}
Let $s\gg 1$ and $t=sN$.
Let $\bC$ be a random configuration of $t$ pebbles on the vertices of $\P_n^d$,
and let $p=(1+\e)s\ln N$ for some $\e>0$.
Then
$$\Pr[\bC(v)<p\ {\rm for\ all\ } v]\ra 1\ {\rm as\ } n\ra\inf\ .$$
\end{lem}

\noindent
{\it Proof.}
Let $q$ be the probability that the vertex $v$ satisfies $\bC(v)\ge p$.
Then $q$ is at most
\begin{eqnarray*}
\frac{{N\pebb t-p}}{{N\pebb t}}
	& = & \frac{t^{\underline p}}{(N+t-1)^{\underline p}}\\
&&\\
& < & \(\frac{t}{N+t}\)^p\\
&&\\
& = & \(1-\frac{1}{s+1}\)^p\\
&&\\
& \le & e^{-p/(s+1)}\ .
\end{eqnarray*}
Hence, the probability that some vertex $v$ satisfies $\bC(v)\ge p$ 
is at most
$$Ne^{-p/(s+1)}\ =\ e^{\ln N(1-\e s)/(s+1)}\ \sim\ N^{-\e}\ \ra\ 0$$
as $n\ra\inf$.
Therefore, the probability that every vertex $v$ satisfies $\bC(v)<p$
tends to 1 as $n\ra\inf$.
\hfill$\pf$
\vspace{0.2 in}

Next we show that a configuration with relatively few pebbles
almost always has some large hole with no pebbles in it.

\begin{lem}\label{hole}
Let $c<2^{-d/(d+1)}$, $u=c(\lg N)^{1/(d+1)}$, $s=2^u$ and $t=sN$.
Write $c=((1-\e)/(2+\d)^d)^{1/(d+1)}$ for some $\e,\d>0$, and set
$m=(2+\d)u$, $M=m^d$ and $k=N/M$.
Partition the vertices of $\P_n^d$ into $k$ disjoint, contiguous blocks
$B_1,\ldots,B_k$ having every side of length $m$.
Let $\bC$ be a random configuration of $t$ pebbles on the vertices of $\P_n^d$.
Then
$$\Pr[\bC(B_h)=0\ {\rm for\ some\ } h]\ra 1\ {\rm as\ } n\ra\inf\ .$$
\end{lem}

\noindent
{\it Proof.}
The second moment method applies.
Let $X_h$ be the indicator variable for the event that the block $B_h$ 
contains no pebbles, and let $X=\sum_{h=1}^kX_h$.
Then Chebyschev's inequality yields
$$\Pr[X=0]\ \le\ \frac{\var[X]}{\bE[X]^2}\ ,$$
and
\begin{eqnarray*}
\var[X] & = & \bE[X^2]-\bE[X]^2 \\
&&\\
& = & \sum_{h,j}\bE[X_hX_j]-\sum_{h,j}\bE[X_h]\bE[X_j] \\
&&\\
& \le & \sum_h\bE[X_h^2]\ ,
\end{eqnarray*}
since $\bE[X_hX_j]\le\bE[X_h]\bE[X_j]$ for $h\neq j$.
Hence,
$$\var[X]\ \le\ \sum_h\bE[X_h^2]\ =\ \sum_h\bE[X_h]\ =\ \bE[X]\ .$$
Moreover, for some $\d>0$ we have
\begin{eqnarray*}
\bE[X] & = & \(\frac{N}{M}\){N-M\pebb t}\Bigg/{N\pebb t}\\
&&\\
& = & \(\frac{N}{M}\)N^{\underline M}\big/(N+t-1)^{\underline M}\\
&&\\
& > & \(\frac{N}{M}\)\(\frac{N-M}{N+t-M}\)^M\\
&&\\
& > & \(\frac{N}{M}\)\(\frac{N-\d N}{N+t-\d N}\)^M\\
&&\\
& = & \(\frac{N}{M}\)\(\frac{1-\d}{1+s-\d}\)^M\\
&&\\
& \sim & \frac{(1-\d)^MN}{Ms^M}\\
&&\\
& = & \frac{(1-\d)^MN^\e}{M}\\
&&\\
& \ra & \inf\ .
\end{eqnarray*}
Hence $\Pr[X=0]\le\var[X]/\bE[X]^2\le 1/\bE[X]\ra 0$ as $n\ra\inf$.
\hfill$\pf$
\vspace{0.2 in}

The following lemma records the structure of the $d$-dimensional grid 
in order to keep track of the results of pebbling steps.

\begin{lem}\label{count}
Let $\bB_m$ be the set of points in $\bZ^d$ whose coordinates are at most
$m/2$ in absolute value, and denote its boundary, those points of $\bB_m$
having some coordinate of absolute value $m/2$, by $\obB_m$.
Define $R_i$ to be the number of points in $\bZ^d-\bB_m$ having distance
$i$ from $\obB_m$, where distance between points is measured by the sum
of absolute values of distances in coordinates.
Then
$$R_i\le \sum_{j=1}^d {d\choose j}2^jm^{d-j}{j\pebb i}\ .$$
\end{lem}

\noindent
{\it Proof.}
We partition the set of points in $\bZ^d-\bB_m$ according to the number
$j$ of coordinates that a given point $\bv$ differs from its nearest
neighbor $\obv$ on $\obB_m$.
Given a fixed $j$, there are $d\choose j$ ways to pick which $j$
coordinates to change.
There are $2^j$ faces of $\obB_m$ to change; positive and negative
coordinates on opposite sides of the origin for each coordinate.
In each case, the coordinates left unchanged on $\obB_m$ can be fixed
at any one of the $m^{d-j}$ values on the given face.
Finally, the number of points at distance $i$ from one of these chosen
pionts on $\obB_m$ equals the number of nonnegative vectors of weight $i$
on j coordinates, namely $j\pebb i$.
(The inequality arises from some slight overcounting due to extra zeroes
that might appear in the distance vectors).
\hfill$\pf$
\vspace{0.2 in}

Finally, our proof of Theorem \ref{grid} in the case of the lower bound
will use this technical lemma to bound the number of pebbles that can
reach the empty hole.

\begin{lem}\label{sum}
$\sum_{i=0}^n {j\pebb i}2^{-i}<2^j\ .$
\end{lem}

\noindent
{\it Proof.}
It is straightforward to use generating functions or induction to prove
$\sum_{i=0}^\inf {j+i-1\choose i}2^{-i}=2^j\ .$
\hfill$\pf$
\vspace{0.2 in}

Turning to the case of the upper bound, we show that almost every
configuration with relatively many pebbles fills every reasonably
large block with plenty of pebbles.

\begin{lem}\label{full}
Let $c^\'=d+1+\e$, some $\e>0$,
$u^\'=c^\'(\lg N)^{1/(d+1)}$, $s^\'=2^{u^\'}$, $t^\'=s^\'N$,
$m^\'=(\frac{\e+1}{c^\'})^{1/d}(\lg N)^{1/(d+1)}$, and $k^\'=N/M^\'$, 
where $M^\'=(m^\')^d$.
Partition the vertices of $\P_n^d$ into $k^\'$ disjoint, contiguous blocks
$B^\'_1,\ldots,B^\'_{k^\'}$ having every side of length $m^\'$.
Let $\bC$ be a random configuration of $t^\'$ pebbles on the vertices 
of $\P_n^d$.
Then
$$\Pr[\bC(B^\'_f)\ge M^\'2^{dm^\'}\ {\rm for\ all\ }f]\
	\ra 1\ {\rm as\ }n\ra\inf\ .$$
\end{lem}

\noindent
{\it Proof.}
We will make use of the fact (see \cite{CEHK}) that every graph $G$ on $V$
vertices has pebbling number at most $V2^\diam(G)$.
Define $Z_f$ to be the event that block $B^\'_f$ contains less than
$M^*=M^\'2^{dm^\'}$ pebbles and approximate the probability
$$\Pr[\cup_{f=1}^k Z_f]\ \le\
	k\sum_{f=0}^{M^*-1}
		{M^\'\pebb f}{N-M^\'\pebb t^\'-f}\Bigg/{N\pebb t^\'}\ .$$
Now use the estimate
$${N-M^\'\pebb t^\'-f}\ \le\ \(\frac{N}{N+t^\'}\)^{M^\'}{N\pebb t^\'}$$
to obtain
$$\Pr[\cup Z_f]\ \le\
        k\(\frac{N}{N+t^\'}\)^{M^\'}\sum_{f=0}^{M^*-1}{M^\'\pebb f}\ .$$
Then use the upper bound 
$$\sum_{f=0}^{M^*-1}{M^\'\pebb f}\ =\
\sum_{f=0}^{M^*-1}{f+1\pebb M^\'-1}\ =\
\sum_{j=1}^{M^*}{j\pebb M^\'-1}\ =\ {M^*\pebb M^\'}\ <\ {M^*}^{M^\'}$$
to obtain
\begin{eqnarray*}
\Pr[\cup Z_f] & < & k\(\frac{N}{N+t^\'}\)^{M^\'}{M^*}^{M^\'}\\
&&\\
& < & \frac{N}{M^\'}\(\frac{M^\'2^{dm^\'}}{s}\)^{M^\'}\\
&&\\
& = & \frac{1}{M^\'}2^{\lg N-M^\'(u-\lg M^\'-dm^\')}\\
&&\\
& = & \frac{1}{M^\'}2^{\lg N-(1+\e)\lg N
	+o(\lg N)+d(\frac{1+\e}{c^\'})^{\frac{d+1}{d}}\lg N}\\
&&\\
& = & \frac{1}{M^\'N^{\e-d(\frac{1+\e}{c^\'})^{\frac{d+1}{d}}-o(1)}}\\
&&\\
& \ra & 0\ 
\end{eqnarray*}
for small enough $\e$.
Thus, almost surely, every $f$ satisfies $\bC(B^\'_f)\ge M^\'2^{dm^\'}$.
\hfill$\pf$
\vspace{0.2 in}

\noindent
{\it Proof of Theorem \ref{grid}.}
We begin with the lower bound.
Given $c<2^{-d/(d+1)}$, we write $c=((1-\e)/(2+\d)^d)^{1/(d+1)}$ 
for some $\e,\d>0$, and set $u=c(\lg N)^{1/(d+1)}$, $s=2^u$, $t=sN$,
$m=(2+\d)u$, $M=m^d$ and $k=N/M$.
Partition the vertices of $\P_n^d$ into $k$ disjoint, contiguous blocks
$B_1,\ldots,B_k$ having every side of length $m$.
Let $\bC$ be a random configuration of $t$ pebbles on the vertices of $\P_n^d$.
By Lemma \ref{hole} we know that, almost surely, some block $B_h$ has 
no pebbles on its vertices.
By Lemma \ref{flat} we know that, almost surely, no other vertex has more
that $p$ pebbles on it, where $p=(1+\e)s\ln N$ for some $\e>0$.

Now, any vertex $v$ with $\bC(v)$ pebbles on it can contribute at most
$\bC(v)/2^i$ pebbles to the boundary $\oB_h$ of $B_h$, where $i$ is
the distance from $v$ to $\oB_h$.
Also, the number of vertices of $P_n^d-B_h$ at distance $i$ from $\oB_h$
is at most $R_i$.
Thus, according to Lemmas \ref{count} and \ref{sum},
the number of pebbles that can be amassed on $\oB_h$
via pebbling steps almost surely is less than
\begin{eqnarray*}
\sum_{i=0}^n pR_i/2^i
& \le & \sum_{i=0}^n p\sum_{j=0}^d{d\choose j}2^jm^{d-j}{j\pebb i}2^{-i}\\
&&\\
& \le & p\sum_{j=0}^d {d\choose j}2^jm^{d-j}\sum_{i=0}^n {j\pebb i}2^{-i}\\
&&\\
& < & p\sum_{j=0}^d{d\choose j}4^jm^{d-j}\\
&&\\
& = & p(m+4)^d\\
&&\\
& \ll & 2^{m/2}\ .\\ 
\end{eqnarray*}
The last line holds because the dominant term in $p(m+4)^d$ is $2^u$, 
and we have $m=(2+\d)u$.
Therefore, almost surely, too few vertices are amassed at distance $m/2$
(or greater) to be able to move a single pebble to the center of $B_h$.
This shows that $\t(\P^d)\in\W(sN)$, as required.

Next we prove the upper bound.
Given $c^\'=d+1+\e$ for some $\e>0$, set $u^\'=c^\'(\lg N)^{1/(d+1)}+2$, 
$s^\'=2^{u^\'}$, $t^\'=sN^\'$, 
$m^\'=(\frac{d+1}{c^\'})^{1/d}(\lg N)^{1/(d+1)}$, 
and $k^\'=N/M^\'$, where $M^\'=(m^\')^d$.
Partition the vertices of $\P_n^d$ into $k^\'$ disjoint, contiguous blocks
$B^\'_1,\ldots,B^\'_{k^\'}$ having every side of length $m^\'$.
Let $\bC$ be a random configuration of $t^\'$ pebbles on the vertices 
of $\P_n^d$.
Then Lemma \ref{full} states that, almost surely, every block $B^\'_f$ 
is full with at least $2^{M^\'}$ pebbles.
Since every graph on $V$ vertices is solvable by $2^V$ pebbles
(see \cite{CEHK}), any given vertex $v$ in $P_n^d$ almost surely is solvable
by the pebbles in the block $B_f$ which contains $v$.
This shows that $\t(\P^d)\in O(s^\'N)$, as required.
\hfill$\pf$

%##############################################################################
%##############################################################################
%
%       REMARKS:
%
\section{Remarks}\label{Remarks}

Theorem 7 brings to mind the following problems, which are of interest
in their own right.

\begin{prob}
Improve the lower bound of $\k(G\Box H)\ge\min\{\d(G),\d(H)\}$.
\end{prob}

\begin{prob}
Find conditions that guarantee $\k(G\Box H)=\d(G)+\d(H)$.
\end{prob}

Let $l=l(n)$ and $d=d(n)$ and denote by $\P_l^d$ the sequence of graphs
$(P_{l(n)}^{d(n)})_{n=1}^\inf$, where $P_l^d=(P_l)^d$.
For $l(n)=2$, $\P_l^n=\Q$, which has a threshold asymptotically less than $N$
by Result \ref{Cubes}.
We conjecture that the same result holds for all fixed $l$.

\begin{conj}\label{cubelike}
Let $\P_l$ denote the graph sequence $(\P_l^n)_{n=1}^\inf$.
Then for fixed $l$ we have $\t(\P_l)\in o(N)$.
\end{conj}

In contrast, we have proved that $\t(\P^d)\in\w(N)$ for fixed $d$.
Thus we believe there should be some relationship between two functions
$l=l(n)$ and $d=d(n)$, both of which tend to infinity, for which the
sequence $\P_l^d$ has threshold on the order of $N$.

\begin{prob}\label{ordern}
Denote by $\P^d$ the graph sequence $(\P_n^{d(n)})_{n=1}^\inf$.
Find a function $d=d(n)\ra\inf$ for which $\t(\P^d)=\Th(N)$.
In particular, how does $d$ compare to $n$?
\end{prob}

%##############################################################################
%##############################################################################
%
%       ACKNOWLEDGMENT:
%
\section*{Acknowledgement}
% (the `*' doesn't number the section)
The authors want to thank Professor Tomasz Luczak for his insight and 
suggestions of how to attack the girth problem.

%##############################################################################
%##############################################################################
%
%       BIBLIOGRAPHY:
%
\bibliographystyle{plain}
%  

%##############################################################################
%##############################################################################
%
%       END OF PAPER
%
\end{document}